\documentstyle[12pt]{article}
\normalbaselines
\begin{document}

\title{Extensions of Rokhlin congruence for curves on surfaces}
\author{Grigory Mikhalkin}
\date{}
\maketitle

\begin{abstract}
The subject of this paper is the problem of arrangement of a real nonsingular
algebraic curve on a  real non-singular algebraic surface. This paper contains
new restrictions on this arrangement extending Rokhlin and
Kharlamov-Gudkov-Krakhnov  congruences for curves on surfaces.
\end{abstract}
\section{Introduction}

If we fix a degree of a real nonsingular algebraic surface then in accordance
with Smith theory there is an upper bound for the total ${\bf Z}_2$-Betti
number of this
surface (the same thing applies to curves). The most interesting case according
to D.Hilbert is the case when the upper bound is reached (in this case the
surface is called an M-surface).

Let $A$  be a real nonsingular algebraic projective curve on a real
   nonsingular
algebraic projective surface $B$. If $A$  is of even degree in $B$  then  $A$
divides $B$  into
two parts $B_+$   and $B_-$   (corresponding to areas of $B$  where polynomial
determining  $A$ is positive and negative).

Rokhlin congruence \cite{R} yields a congruence modulo 8 for Euler
characteristic
$\chi$ of $B_+$ provided that \begin{description}
\item[$(i)$]$B$ is an M-surface \item[$(ii)$] $A$ is an M-curve
\item[$(iii)$] $B_+$ lies wholly in one component of $B$ \item[$(iv)$] $
rk(in_*:H_1(B_+;{\bf Z}_2)\rightarrow H_1(B;{\bf Z}_2))=0 $ \item[$(v)$] if
the degree of polynomial determining
$A$ in $B$ is congruent to 2 modulo 4 then all components of $B$ containing no
components of $A$ are contractible in $P^q$ \footnote{{\em\underline {Remark}}.
 Paper \cite{R} contains a mistake in the calculation of characteristic class
of
covering $Y\rightarrow {\bf C}B$. It leads to the omission (after
reformulation) of $(v)$ in assumptions of congruence. The proof given in
\cite{R}
really uses $(v)$.
} \end{description}
Kharlamov-Gudkov-Krakhnov \cite{Kh},\cite{GK} congruence yields congruence
modulo 8
for $\chi(B_+)$ under assumptions $(i)$,$(iii)$,$(v)$ and either assumption
that $A$ is an
(M-1)-curve and $rk(in_*)=0$ or assumption that $A$ is an M-curve and
$rk(in_*)=1$ and all
components of $A$ are ${\bf Z}_2$-homologically trivial. Recall that Rokhlin
congruence for surfaces yields congruence modulo 16 for $\chi(B)$ so a
congruence modulo 8 for $\chi(B_+)$ is equivalent to a congruence modulo 8
for $\chi(B_-)$.

One of the properties of M-surfaces (similar to the property of M-curves)
 remarked by V.I.Arnold \cite{A}
is that a real M-surface is a characteristic surface
 in its own complexification (for M-curves it means that a real M-curve divides
its own complexification since a complex curve is orientable).We shall
say that a real surface is of a characteristic type if it is a characteristic
surface in its own complexification.Note that the notion of characteristic
type of surfaces is analogous to the notion of type I of curves.

Consider at first the weakening of assumption $(i)$ in Rokhlin and
Kharlamov-Gudkov-Krakhnov congruences.Instead of $(i)$ we can only assume
that B is of characteristic type.

According to O.Ya.Viro \cite{V} there are some extra structures on real
surfaces
of characteristic type ,namely, $Pin_-$-structures and semiorientations or
relative semiorientations (semiorientation is the orientation up to the
reversing). In this paper we introduce another structure on surfaces of
characteristic type --- complex separation which is also determined by the
arrangement of a real surface in its complexification. The complex separation
is a natural separation of the set of components of a real surface of
characteristic type into two subsets.Note that the set of semiorientations
of a surface is an affine space over ${\bf Z}_2$-vector space of separations
of this surface.

We use the complex separation to weaken assumption $(iii)$. In theorem 1
instead
of $(iii)$ we assume only that $B_+$ lies in components of one class of
complex separation.

The further extension ,theorem 2, can be applied not only for curves of even
degree but also sometimes for  curves of odd degree. Another weakening of
assumptions in theorem 2 is that components of a curve are not necessarily
${\bf Z}_2$-homologically trivial.

These extensions can be applied to curves on quadrics and cubics.
Theorem 1 together with an analogue of Arnold inequality for curves on cubics
gives a complete system of restrictions for real schemes of flexible curves
of degree 2 on cubics of characteristic type (see \cite{M1}).An application of
theorem 2 to curves on an ellipsoid gives a complete system of restrictions
for real schemes of flexible curves of degree 3 on an ellipsoid and reduces
the problem of classification of real schemes of flexible curves of degree 5
on an ellipsoid to the problem of the existence of two real schemes (see
\cite{M2}).
An application of theorem 2 to curves on a hyperboloid extends Matsuoka
congruences \cite{Ma} for curves with odd branches on a hyperboloid (see
\cite{M3}).

Applications of theorem 2 to empty curves on surfaces give restrictions for
surfaces involving complex separation of surfaces. Restrictions for curves
of even degree on surfaces can be obtained also by the application of these
restrictions for surfaces to the 2-sheeted covering of surface branched along
the curve , if we know the complex separation of this covering. This complex
separation is determined by complex orientation of the curve.For example in
this way one can obtain new congruences for complex orientations of curves on
a hyperboloid.

These applications and further extension of Rokhlin congruence for curves on
surfaces will be published separately in \cite{M3}. For example the assumption
that
the surface and the curve on the surface are complete intersections is quite
unnecessary ,but this assumption simplifies definition of number $c$ in
formulations of theorems. The formulations of these results were announced
in \cite{M2} as well as the formulations of results of the present paper.

The author is indebted to O.Ya.Viro for advices.
\section{Notations and formulations of main theorems}
Let the surface $B$ be the transversal intersection of hypersurface in
$P^q$ defined by equations $P_j(x_0,\ldots,x_q)=0 ,j=1,\ldots,s-1$ ;
${\bf C}B$ and ${\bf R}B$ be sets of complex and real points of $B$;
let $A$ be a nonsingular curve on $B$ defined by an equation
$P_s(x_0,\ldots,x_q)=0$ ,where $P_j $ is a real homogeneous polynomial
of degree $m_j$ $j=1,\ldots,s$, $q=s+1$; ${\bf C}A$ and ${\bf R}A$ be sets of
complex
and real points of $A$. Let $conj$ denote the involution of complex
conjugation.
 Set $c$ to be equal to $\frac{\prod_{j=1}^{s-1} m_{j}}{4}$.If $m_s$ is
even then denote \{$x\in {\bf R}B | \pm P_s(x)\geq 0$\} by $B_{\pm}$ and
set $d$ to be equal to $rk(in_*:H_1(B_+;{\bf Z}_2)\rightarrow H_1({\bf R}B;{\bf
Z}_2))$.

A real algebraic variety is called an (M-$j$)-variety if its total ${\bf
Z}_2$-Betti number is less by $2j$ then total ${\bf Z}_2$-Betti number of
its complexification (Harnack-Smith inequality shows that $j\geq 0$). Let
$A$ be an (M-$k$)-curve. Let $D_M$ be the operator of Poincar\'{e} duality of
manifold $M$.We shall say that $B$ is a surface of characteristic type if
$[{\bf R}B]=D_{{\bf C}B}w_2({\bf C}B)\in H_2({\bf C}B;{\bf Z}_2)$
(as it is usual we denote by $[{\bf R}B]$ the element of
$H_2({\bf C}B;{\bf Z}_2)$ realized by ${\bf R}B$).
We shall say that $(B,A)$ is a pair of characteristic type if
$[{\bf R}B]+[{\bf C}A]+D_{{\bf C}B}(w_2({\bf C}B))=0\in H_2({\bf C}B;{\bf
Z}_2)$.It is said that $A$ is a curve of type I if $[{\bf R}A]=0\in H_1({\bf
C}A;{\bf Z}_2)$.It is said that $A$ is of even(odd) degree if
$[{\bf C}A]=0\in H_2({\bf C}B;{\bf Z}_2)$ (otherwise).

As it is usual we denote by $\sigma$ and $\chi$ the signature and the Euler
characteristic.By $\beta(q)$ we mean Brown invariant of ${\bf Z}_4$-valued
quadratic
form $q$.

\newtheorem{theorem}{Theorem}
\begin{theorem}
If $B$ is a surface of characteristic type then there is defined a natural
separation of surface ${\bf R}B$ into two closed surfaces $B_1$ and $B_2$
by the condition that $B_j,j=1,2$, is a characteristic surface in
${\bf C}B/conj$.Suppose that $m_s$ is even,$B_+\subset B_1$,
every component of ${\bf R}A$ is ${\bf Z}_2$-homologous to zero in ${\bf R}B$
and if $m_s\equiv 2\pmod{4}$ then suppose besides that $B_2$ is contractible
in ${\bf R}P^q$.
\begin{itemize}\begin{description} \item[a)] If $d+k=0$ then $\chi(B_+)\equiv
c\pmod{8}$.
\item[b)] If $d+k=1$ then $\chi(B_+)\equiv c\pm 1\pmod{8}$.
\item[c)] If $d+k=2$ and $\chi(B_+)\equiv c+4\pmod{8}$ then $A$ is of type I
and $B_+$ is orientable.
\item[d)] If $A$ is of type I and $B_+$ is orientable then $\chi(B_+)\equiv
c\pmod{4}$.
\end{description}\end{itemize}\end{theorem}
\begin{theorem}
If $(B,A)$ is a pair of characteristic type then there is defined a natural
separation of surface ${\bf R}B\setminus {\bf R}A$ into two surfaces
$B_1$ and $B_2$ with common boundary ${\bf R}A$ by the condition that
$B_j\cup{\bf C}A/conj$ is a characteristic surface in ${\bf C}B/conj$,
there is defined Guillou-Marin form $q_j$ on $H_1(B_j\cup{\bf C}A/conj;{\bf
Z}_2)$ and \begin{displaymath}
\chi(B_j)\equiv c+\frac{\chi({\bf R}B)-\sigma({\bf C}B)}{4}+\beta(q_j)\pmod{8}
   .\end{displaymath}\end{theorem}
\section{Proof of theorem 2}
Consider the Smith exact sequence of double branched covering
$\pi:{\bf C}B\rightarrow {\bf C}B/conj$\begin{displaymath}
\stackrel{\beta_3}{\rightarrow}H_3({\bf C}B/conj,{\bf R}B;{\bf
Z}_2)\stackrel{\gamma_3}{\rightarrow}H_2({\bf R}B;{\bf Z}_2)\oplus H_2({\bf
C}B/conj,{\bf R}B;{\bf Z}_2)\stackrel{\alpha_2}{\rightarrow}H_2({\bf C}B;{\bf
Z}_2)\stackrel{\beta_2}{\rightarrow}  .\end{displaymath}
Let $\phi$ denote the composite homomorphism \begin{displaymath} H_2({\bf
C}B/conj,{\bf R}B;{\bf Z}_2)\stackrel{0\oplus id}{\rightarrow}H_2({\bf R}B;{\bf
Z}_2)\oplus H_2({\bf C}B/conj,{\bf R}B;{\bf
Z}_2)\stackrel{\alpha_2}{\rightarrow}H_2({\bf C}B;{\bf Z}_2)
.\end{displaymath}
Let $j$ denote the inclusion map $({\bf C}B/conj,\emptyset)\rightarrow({\bf
C}B/conj,{\bf R}B)$.Recall that $\phi_*\circ j_*$ is equal to Hopf homomorphism
$\pi^!$. It is easy to deduce from the exactness of the Smith sequence that
$\phi$ is a monomorphism. Indeed, $\pi_1({\bf C}B)=0$ hence $\pi_1({\bf
C}B/conj)=0$ and $H_3({\bf C}B/conj;{\bf Z}_2)=0$.
Therefore boundary homomorphism $\partial:H_3({\bf C}B/conj,{\bf R}B;{\bf
Z}_2)\rightarrow H_2({\bf R}B;{\bf Z}_2)$ is a monomorphism.Therefore,since
$\partial$
is the first component of $\gamma_3$, $Im \gamma_3\cap(\{0\}\oplus H_2({\bf
C}B/conj,{\bf R}B;{\bf Z}_2))=0$ and $\phi$ is a monomorphism.

It is easy to check that \begin{displaymath}\pi^*w_2({\bf C}B/conj)=D_{{\bf
C}B}^{-1}[{\bf R}A]+w_2({\bf C}B)  .\end{displaymath} Thus $\pi^!(D_{{\bf
C}B/conj}w_2({\bf C}B/conj))=[{\bf C}A]$,
therefore,because of the injectivity of $\phi$,we obtain that
\begin{displaymath} j_*D_{{\bf C}B/conj}w_2({\bf C}B/conj)=[{\bf C}A/conj,{\bf
R}A]\in H_2({\bf C}B/conj,{\bf R}B;{\bf Z}_2)  .\end{displaymath} It means that
there exists a compact surface
$B_1\subset{\bf R}B$ with boundary ${\bf R}A$ such that
$B_1\cup{\bf C}A/conj$ is a characteristic surface in ${\bf C}B/conj$.
Surface ${\bf R}A$ is homologous to zero in ${\bf C}B/conj$
since ${\bf R}A$ is the set of branch points of $\pi$.
Set $B_2$ to be equal to the closure of $({\bf R}B\setminus B_1)$.
Then $B_2\cup {\bf C}A/conj$ is a characteristic surface ${\bf C}B/conj$,
$B_1\cup B_2={\bf R}B$,$B_1\cap B_2=\partial B_1=\partial B_2={\bf R}A$.

Let us prove the uniqueness of pair $\{B_1,B_2\}$.It is sufficient to prove
that the dimension of the kernel of inclusion homomorphism
$H_2({\bf R}B;{\bf Z}_2)\rightarrow H_2({\bf C}B/conj;{\bf Z}_2)$
is equal to 1. This follows from the equality
$\dim H_3({\bf C}B/conj,{\bf R}B;{\bf Z}_2)=1$ that can be deduced from
the exactness of the Smith sequence.

We apply now Guillou-Marin congruence \cite{GM} to pair $({\bf
C}B/conj,B_j\cup{\bf C}A/conj),j=1,2$\begin{displaymath}
\sigma({\bf C}B/conj)\equiv[B_j\cup{\bf C}A/conj]\circ[B_j\cup{\bf
C}A/conj]+2\beta(q_j)\pmod{16}  .\end{displaymath}
Hirzebruch index theorem gives an equality $\sigma({\bf
C}B/conj)=\frac{\sigma({\bf C}B)-\chi({\bf R}B)}{2}$.To finish the proof note
that
$[B_j\cup{\bf C}A/conj]\circ[B_j\cup{\bf C}A/conj]=2c-2\chi(B_j)$
(the calculation is similar to Marin calculation in \cite{Marin}).
\section{Proof of the theorem 1}
Pair $(B,\emptyset)$ is of characteristic type since $A$ is of even degree
in $B$. Thus the first part of theorem 1 follows from theorem 2 --- there exist
a natural separation of $B$ into two surfaces $B_1$ and $B_2$ such that
$B_1$ and $B_2$ are characteristic surfaces in ${\bf C}B/conj$.

Let $V$ denote $B_+\cup{\bf C}A/conj$.Let $W$ denote $V\cup B_1$.
Recall that $B_+\cap B_2=\emptyset$ thus $V\cap B_2=\emptyset$.
\newtheorem{lemma}{Lemma}
\begin{lemma}
$[V]=0\in H_2({\bf C}B/conj;{\bf Z}_2)$\end{lemma}
{\em\underline{Proof}}. Since $A$ is of even degree in $B$, there exists a
2-sheeted covering $p:Y\rightarrow{\bf C}B$ branched along ${\bf C}A$.
Involution $conj$ can be lifted to involutions $T_+$ and $T_- :Y\rightarrow Y$
since ${\bf C}A$ is invariant under $conj$. It is easy to see using the
straight
algebraic construction of $p$ that $T_+$ and $T_-$ can be chosen in such
a way that the set of fixed points of $T_\pm$ is $p^{-1}(B_\pm)$.

Consider the diagram\begin{picture}(200,-100)(-50,0)
\put(0,0){$Y$}\put(60,0){${\bf C}B$}
\put(0,-40){$Y/T_-$}\put(60,-40){${\bf C}B/conj .$}
\put(10,4){\vector(1,0){46}}
\put(3,-3){\vector(0,-1){25}}
\put(67,-3){\vector(0,-1){25}}
\put(33,10){$p$}\put(71,-15){$\pi$}
\end{picture}\\[42pt]
This diagram can be expanded to a commutative one by map
$p^{'}:Y/T_\mp\rightarrow{\bf C}B/conj$.
It is easy to see that $p^{'}$ is a 2-sheeted covering branched along $V$.
Therefore $[V]=0\in H_2({\bf C}B/conj;{\bf Z}_2)$.

Using Lemma 1 we see that $W$ is a characteristic surface in ${\bf C}B/conj$
as well as $B_2$.We apply Guillou-Marin congruence to these two surfaces :
\begin{displaymath}
\sigma({\bf C}B/conj)\equiv[W]\circ[W]+2\beta(q_W)\equiv
2c-2\chi(B_+)-2\chi(B_2)+2\beta(q_W)\pmod{16}\end{displaymath}
\begin{displaymath}\sigma({\bf
C}B/conj)\equiv[B_2]\circ[B_2]+2\beta(q_{B_2})\equiv
-2\chi(B_2)+2\beta(q_{B_2})\pmod{16}\end{displaymath}
,where $q_W$ and $q_{B_2}$ are Guillou-Marin forms of $W$ and $B_2$.
Therefore \begin{displaymath} \chi(B_+)\equiv
c+\beta(q_W)-\beta(q_{B_2})\pmod{8}  .\end{displaymath}
\begin{lemma}
$\forall x\in H_1(B_2;{\bf Z}_2), q_{B_2}(x)-q_W(x)=\left\{ \begin{array}{ll}
0 & \mbox{if $x$ is contractible in ${\bf R}P^q$} \\
\frac{m_s}{2} & \mbox{if $x$ is noncontractible in ${\bf R}P^q$ .}
\end{array} \right. $
\end{lemma}
{\em\underline{Proof}}. It follows from the definition of Guillou-Marin
form that values on $x$ of $q_{B_2}$ and $q_W$ are differed by linking
number of $x$ and $V$ in ${\bf C}B/conj$ that is equal to linking
number of $x$ and ${\bf C}A$ in ${\bf C}B$. The last linking number can be
calculated from the straight construction of a 2-sheeted covering
branched along ${\bf C}A$.

It was shown in \cite{KV} that Brown invariant of form $q$ on the union of two
surfaces with common boundary is equal to the sum of Brown invariants
of restrictions of $q$ on these surfaces in the case when $q$ vanishes on
the common boundary.
Now theorem 1 follows from this additivity of Brown invariant
and the classification of low-dimensional ${\bf Z}_4$-valued quadratic forms
(see\cite{KV}). Indeed, since every component of ${\bf R}A$ is homologous to 0
in ${\bf R}B$,
$\beta(q_W)=\beta(q_W{|_{{\bf C}A/conj}})+\beta(q_W|_{B_+})+\beta(q_W|_{B_2})$.
Lemma 2 shows that under assumptions of theorem 1
$\beta(q_W|_{B_2})=\beta(q_{B_2})$.
To complete the proof note that ranks of intersection forms on
$H_1(B_+;{\bf Z}_2)$ and $H_1({\bf C}A/conj;{\bf Z}_2)$ are equal to
$d$ and $k$ respectively.

\end{document}